\documentclass[10pt]{article}
\usepackage[a4paper,hmargin=27mm,vmargin=22mm]{geometry}
\usepackage{amsmath,amssymb,amsthm,mathtools,microtype}
\usepackage[colorlinks=true,linkcolor=blue,citecolor=blue,urlcolor=blue]{hyperref}
\hypersetup{pdfauthor={Volkov Sergey},
 pdftitle={A strict upper bound for three-point weighted selection in planar vector-valued regression}}
\newtheorem{theorem}{Theorem}
\newtheorem{lemma}[theorem]{Lemma}
\newtheorem{proposition}[theorem]{Proposition}
\newcommand{\R}{\mathbb R}
\newcommand{\tr}{\operatorname{tr}}
\newcommand{\adj}{\operatorname{adj}}
\newcommand{\conv}{\operatorname{conv}}
\newcommand{\supp}{\operatorname{supp}}
\newcommand{\Fw}{F_{\mathrm{weighted}}(2,2,3)}
\newcommand{\T}{\mathcal T_3}
\title{A strict upper bound for three-point weighted selection\\in planar vector-valued regression}
\author{Volkov Sergey}
\begin{document}
\maketitle
\begin{abstract}
We prove $F_{\mathrm{weighted}}(2,2,3)<11/6$ for minimum-Frobenius-norm
least-squares regression with two input and two output coordinates,
improving the upper bound $15/8$. More precisely, a universal
$\varepsilon_0>0$ satisfies
$\mathcal T_3\le(5/6-\varepsilon_0)E$ for every finite residual frame.
The proof combines an exact scalar certificate, rigidity of equality in
planar selection, and a compactness argument that controls disappearing
weights. Together with the known lower bound, this gives
$13/8\le F_{\mathrm{weighted}}(2,2,3)<11/6$.
No numerical value of $\varepsilon_0$ is claimed.
\end{abstract}

\section{Problem and result}
Weighted data selection for regression was proposed by
Hanneke, Moran, Shlimovich and Yehudayoff~\cite{open,erms}.
Zhang~\cite{zhang} established
$13/8\le\Fw\le15/8$ and reduced the conjectured value $13/8$
to a finite planar moment problem. We use Zhang's residual-frame formulation,
interpolation hulls, star localization and seven-atom reduction.
An exact scalar certificate, rigidity of equality and compactness yield
a strict uniform improvement over all finite datasets.
The conjectured value $13/8$ remains open.

For a finite sequence $D=((x_i,y_i))_{i=1}^N$ in $\R^2\times\R^2$, $N\ge1$, let
$L_D(W)=N^{-1}\sum_i\|Wx_i-y_i\|^2$ and $L_D^*=\min_W L_D(W)$.
For nonnegative weights $c$ of sum one, let $\widehat W_c$ be the
minimum-Frobenius-norm minimizer of $\sum_i c_i\|Wx_i-y_i\|^2$.
Then $\Fw$ is the supremum over datasets of
\[
 \frac{\inf_{c_i\ge0,\,\sum_i c_i=1,\,|\supp c|\le3}L_D(\widehat W_c)}{L_D^*},
\]
with $0/0=1$ and positive-over-zero equal to $+\infty$.
We take the infimum because the selected fit can jump when the weighted
design loses rank. A limit of full-rank fits bounds the infimum via the
continuity of $L_D$ as a function of $W$.

A residual frame is a finite collection $p_i>0$, $u_i,w_i\in\R^2$ with
\begin{equation}\label{frame}
 \sum_i p_i=1,\quad\|u_i\|^2=2,\quad
 \sum_i p_i u_i u_i^\top=I_2,\quad\sum_i p_i w_i u_i^\top=0.
\end{equation}
Put $E=\sum_i p_i\|w_i\|^2$ and $R=\sum_i p_iw_iw_i^\top$.
For independent $u_i,u_j$, let $G_{ij}u_i=w_i$, $G_{ij}u_j=w_j$.
Using Frobenius distance, define
\[
 \T=\min_{\substack{2\le|S|\le3\\S\text{ contains an independent pair}}}
 \operatorname{dist}\left(0,
 \conv\{G_{ij}:i,j\in S,\ \det(u_i,u_j)\ne0\}\right)^2.
\]
Tightness guarantees an independent pair.
\begin{theorem}\label{main}
There exists a universal $\varepsilon_0>0$ such that every residual frame
satisfies $\T\le(5/6-\varepsilon_0)E$. Consequently
$13/8\le\Fw\le11/6-\varepsilon_0<11/6$.
\end{theorem}

\section{Planar selection and the scalar certificate}
\begin{lemma}\label{planar}
Let $v_j\in\R^2$ have positive probabilities $p_j$ and mean zero. For any
real $2\times2$ matrix $A$ with $\tr A=1$, some convex combination $v$ of
at most two of these points satisfies
\[
 \|v\|^2\le\sum_jp_j\|Av_j\|^2.
\]
\end{lemma}
\begin{proof}
Minimize the right-hand linear objective over centered probability vectors
on the fixed points, and choose a minimizer $q$ of least support. If the
columns $(1,v_j)$ on its positive support are dependent, a supported kernel
vector has both signs. Small perturbations in either direction remain
feasible, so optimality forces zero objective slope. Continuing until one
weight vanishes contradicts least support. Hence at most three affinely
independent points remain and their objective is no larger than the original.
If at most two remain, their centered mean itself is zero.

Otherwise let $S_q=\sum_jq_jv_jv_j^\top$. It is positive definite: a null
direction would place all three points on a line. Put $n_j=S_q^{-1}v_j$,
choose $i$ maximizing $r=\|n_i\|^2>0$, and set
\[
 \lambda_j=q_j\left(1-\frac{n_j^\top n_i}{r}\right).
\]
Cauchy--Schwarz and maximality give $\lambda_j\ge0$. Centering gives
$\sum_j\lambda_j=1$, and $\lambda_i=0$. Since
$\sum_jq_jv_jn_j^\top=I_2$,
\[
 v=\sum_j\lambda_jv_j=-\frac{n_i}{r},\qquad \|v\|^2=\frac1r.
\]
Finally, writing $C_q=\sum_jq_j\|Av_j\|^2$,
\[
 1=\sum_jq_j n_j^\top Av_j,\qquad
 1\le\left(\sum_jq_j\|n_j\|^2\right)C_q\le rC_q.
\]
Thus $\|v\|^2\le C_q$, proving the assertion.
\end{proof}

For positive semidefinite $S$ with $\tr S>0$, the matrix
$A=\adj S/\tr S$ has trace one, and in dimension two
\begin{equation}\label{adjidentity}
 \tr(ASA^\top)=\frac{\det S}{\tr S}.
\end{equation}
This identity holds also when $S$ is singular. If a covariance has trace
zero, all its positive-weight points vanish.

Fix an anchor $i$ and a unit vector $n_i\perp u_i$. Put
\[
 M_i=\tfrac12w_iu_i^\top,\quad
 b_{ij}=u_j^\top n_i,\quad c_{ij}=\tfrac12u_i^\top u_j,
 \quad \xi_{ij}=\frac{w_j-c_{ij}w_i}{b_{ij}}\quad(b_{ij}\ne0).
\]
Because $u_j=c_{ij}u_i+b_{ij}n_i$,
\[
 G_{ij}=M_i+\xi_{ij}n_i^\top,\qquad
 \|M_i+\xi n_i^\top\|_F^2=\tfrac12\|w_i\|^2+\|\xi\|^2.
\]
The cross term vanishes since $u_i^\top n_i=0$. The positive star weights
$P_{ij}=p_jb_{ij}^2$ satisfy
\[
 \sum_{b_{ij}\ne0}P_{ij}=1,\qquad
 \sum_{b_{ij}\ne0}P_{ij}\xi_{ij}=0.
\]
Indeed these are respectively $n_i^\top I_2n_i=1$ and the difference
between the cross-moment $\sum_jp_jw_ju_j^\top n_i=0$ and
$w_i u_i^\top I_2n_i/2=0$.

Define
\[
 L_i=\sum_{b_{ij}=0}p_j(w_j-c_{ij}w_i)(w_j-c_{ij}w_i)^\top\succeq0,
 \qquad H_i=R+\tfrac12w_iw_i^\top.
\]
The star covariance is
\begin{equation}\label{starcov}
 S_i=\sum_{b_{ij}\ne0}P_{ij}\xi_{ij}\xi_{ij}^\top=H_i-L_i.
\end{equation}
To check it, expand the same sum over all $j$. The two cross terms vanish
because $\sum_jp_jc_{ij}w_j=0$, and $\sum_jp_jc_{ij}^2=1/2$.
Subtracting the parallel terms yields \eqref{starcov}.

Apply Lemma~\ref{planar} with $A_i=\adj H_i/\tr H_i$; its denominator
$E+\|w_i\|^2/2$ is positive for $E>0$. PSD monotonicity and
\eqref{adjidentity} give a two-neighbor convex mean with
\[
 \|\xi\|^2\le\tr(A_iS_iA_i^\top)
 =\frac{\det H_i}{\tr H_i}-\tr(A_iL_iA_i^\top).
\]
Consequently, for every anchor,
\begin{equation}\label{anchor}
 \T\le \frac{\|w_i\|^2}{2}
 +\frac{\det(R+w_iw_i^\top/2)}{E+\|w_i\|^2/2}
 -\tr(A_iL_iA_i^\top).
\end{equation}
This involves at most three atoms. Dropping the last, nonnegative term is
valid but potentially loses information.

Normalize $E=1$ and rotate output coordinates so $R=\operatorname{diag}(r,1-r)$.
Set $a_i=w_{i1}^2$, $b_i=w_{i2}^2$, and $s=a+b$. The first two terms of
\eqref{anchor} are
\[
 g_r(a,b)=\frac{s}{2}+
 \frac{r(1-r)+((1-r)a+rb)/2}{1+s/2},\qquad
 \sum_i p_ia_i=r,\quad\sum_i p_ib_i=1-r.
\]
Put $\delta=2r-1$, $x=1+3\delta/10$, $y=1-3\delta/10$,
$P=xy$, $Q=7/3-21\delta^2/100$, and $h=P-ay-bx$.
Since
\[
 \sum_i p_i h(a_i,b_i)=\frac{21}{100}\delta^2\ge0,
\]
some positive-weight atom has $h\ge0$. Define
\[
 F(a,b)=s^2+\frac43s-\delta(a-b)-\frac73-\delta^2.
\]
Exact expansion gives
\begin{align}
 PF&=-\frac{121}{100}\delta^2P
     -h(Q+ax+by)-ab(x-y)^2,\label{poly}\\
 4(1+s/2)(g_r(a,b)-5/6)&=F(a,b).\label{gidentity}
\end{align}
Here $x,y,P,Q>0$ for $|\delta|\le1$. At the selected atom all terms on the
right of \eqref{poly} are nonpositive. This proves $\T\le5/6$.
Homogeneity gives $\T\le5E/6$ for general energy.

We also need its quantitative anisotropy consequence. At $h\ge0$,
$s\le P/\min(x,y)=\max(x,y)\le13/10$. Thus \eqref{poly}--\eqref{gidentity}
give
\begin{equation}\label{anisogap}
 \T\le\frac56-\frac{11}{60}\delta^2.
\end{equation}
The coefficient is $(121/100)/(4(1+13/20))=11/60$.

The scalar relaxation alone is sharp: the one-atom scalar distribution
$r=a=b=1/2$ has $g=5/6$. It is not a valid tight frame.
The strict improvement must therefore use geometry.

\section{Equality cannot occur}
\begin{lemma}[Equality rigidity]\label{rigidity}
Suppose the centered positive distribution has positive definite covariance
$S$, put $C=\det S/\tr S>0$ and $A=\adj S/\tr S=CS^{-1}$, and let $d^2$
be the minimum squared norm of a convex combination of at most two original
points. If $d^2=C$, then every original point satisfies
\[
 \|Av_j\|^2=C,\qquad \|S^{-1}v_j\|^2=\tr(S^{-1})=1/C.
\]
\end{lemma}
\begin{proof}
The original $A$-objective is $C$ by \eqref{adjidentity}. Minimize it as in
Lemma~\ref{planar}. Its minimum $C_q$ obeys $d^2\le C_q\le C$, so all
equal $C$. A minimizing least-support distribution must be a centered
positive triangle, since otherwise $d=0$.

In its planar proof all inequalities are equalities:
\[
 C=d^2\le\frac1r\le
 \frac1{\sum_jq_j\|S_q^{-1}v_j\|^2}\le C_q=C.
\]
Equality in the maximum step forces all three dual norms to equal $r=1/C$.
Equality in weighted Cauchy--Schwarz forces
$Av_j=C S_q^{-1}v_j$ on the triangle. Its vertices span $\R^2$, hence
$A=C S_q^{-1}$, and comparison with $A=CS^{-1}$ gives $S_q=S$.
Each of its three costs $\|Av_j\|^2$ therefore equals $C$.

For an original point $v_k$ not in the selected support, write its affine
coordinates in the triangle (which need not be nonnegative):
$v_k=\sum_{j=1}^3\beta_jv_j$, $\sum_j\beta_j=1$. For sufficiently small
$t>0$, put weight $t$ on $v_k$ and change triangle weights to
$q_j-t\beta_j\ge0$. This is still a centered probability vector. Its cost
is $C+t(\|Av_k\|^2-C)$. Minimality forces $\|Av_k\|^2\ge C$.
Every original weight is positive and their costs average $C$, so every
cost is exactly $C$. The second identity follows from $A=CS^{-1}$ and
$\tr S^{-1}=\tr S/\det S$.
\end{proof}

\begin{lemma}\label{angles}
A tight planar frame contains two indices with
$c^2=(u_i^\top u_j/2)^2\le1/4$.
\end{lemma}
\begin{proof}
Write $u_i/\sqrt2=(\cos\theta_i,\sin\theta_i)$. Tightness implies that
the unit-circle points $t_i=(\cos2\theta_i,\sin2\theta_i)$ have weighted
mean zero. Reduce a convex representation of zero to at most three points.
If two suffice, they are antipodal and the corresponding $c^2=0$.
Otherwise their triangle contains zero in its interior; its three cyclic
angular gaps are at most $\pi$ and sum $2\pi$. Some gap lies in
$[2\pi/3,\pi]$, so that pair has
$c^2=(1+\cos(2\theta_i-2\theta_j))/2\le1/4$.
A boundary triangle likewise has an antipodal pair.
\end{proof}

\begin{proposition}\label{strictfinite}
Every finite positive-energy residual frame satisfies $\T<5E/6$.
\end{proposition}
\begin{proof}
Normalize $E=1$ and suppose $\T=5/6$. Equation~\eqref{anisogap} forces
$R=I_2/2$. In this case the dropped-correction anchor expression is
\[
 f(s)=s/2+(1+s)/(4+2s),\qquad f(1)=5/6,\quad
 f'(s)=1/2+1/(2(s+2)^2)>0.
\]
Since $\T\le f(\|w_i\|^2)$ for every $i$, all $\|w_i\|^2\ge1$.
Their weighted mean is one and weights are positive, hence all equal one.

For every anchor, $H_i=I_2/2+w_iw_i^\top/2$ has eigenvalues $1,1/2$,
so $A_i$ is invertible and $\det H_i/\tr H_i=1/3$.
Equation~\eqref{anchor} forces $\tr(A_iL_iA_i^\top)=0$, hence $L_i=0$
and $S_i=H_i$. The minimum squared norm of a two-point convex mean in this
star must equal $1/3$: Lemma~\ref{planar} bounds it above by $1/3$,
while any smaller value would produce a certificate below $5/6$.
Lemma~\ref{rigidity} therefore forces, for every star neighbor,
\begin{equation}\label{dualnorm}
 \|H_i^{-1}\xi_{ij}\|^2=\tr(H_i^{-1})=3.
\end{equation}

Choose $i,j$ from Lemma~\ref{angles}, and set
$c=u_i^\top u_j/2$, $d=w_i^\top w_j$. This is an independent pair.
In output coordinates with first axis $w_i$, the two components of
$w_j-cw_i$ are $d-c$ and a number of squared magnitude $1-d^2$.
Also $b_{ij}^2=2(1-c^2)$, and $H_i^{-1}$ has eigenvalues $1,2$.
Thus \eqref{dualnorm} is equivalent to
\[
 \frac{(d-c)^2+4(1-d^2)}{2(1-c^2)}=3,\qquad
 7c^2-2cd-3d^2=2.
\]
But exact square completion gives
\[
 7c^2-2cd-3d^2=\frac{22}{3}c^2-3(d+c/3)^2
 \le\frac{11}{6}<2,
\]
a contradiction. This proves pointwise strictness.
\end{proof}

\section{A uniform gap}
\begin{lemma}\label{seven}
For each positive-energy frame there is a residual frame on a subset of
at most seven of its fixed atoms, with energy $E'\le E$, and with
$\T'\ge\T$. If $E'=0$, then $\T=0$.
\end{lemma}
\begin{proof}
We recall Zhang's support reduction~\cite{zhang}.
Minimize $\sum_i q_i\|w_i\|^2$ over nonnegative weights on the fixed
atoms satisfying the frame constraints. The symmetric input covariance gives three linear equations and the
cross-covariance gives four. Normalization follows by taking the trace
of the input equation, since $\|u_i\|^2=2$. Thus there are at most
seven independent linear constraints.
The feasible probability polytope is compact and nonempty. At a minimizing
weight vector of least support, a dependence of its active constraint
columns would permit both signs of small feasible perturbation. The energy
slope must be zero by optimality. Continuing to the first zero weight would
reduce support while preserving energy, a contradiction. The positive
support consequently has at most seven atoms.
The pair interpolators do not depend on the frame weights. Deleting atoms
only removes candidates from the selection minimum, giving $\T'\ge\T$.
If $E'=0$, all remaining $w_i$ vanish; tightness still supplies an independent
pair and its interpolator is zero. Then $\T' =0$, so $\T=0$.
\end{proof}

\begin{proof}[Proof of the uniform frame assertion in Theorem~\ref{main}]
Let $\Gamma=\sup\T/E$ over all finite positive-energy frames.
We know $\Gamma\le5/6$. Suppose $\Gamma=5/6$.
For each $k$ choose a frame with positive ratio at least $5/6-\epsilon_k$,
where $\epsilon_k\downarrow0$. Lemma~\ref{seven} gives a reduced frame with
$E'>0$ and $\T'/E'\ge\T/E$. Scale the reduced $w_i$ by $1/\sqrt{E'}$.
The resulting frames have energy one and $\T^{(k)}\ge5/6-\epsilon_k$.
There are at most seven atoms. Relabel and pass to a subsequence of constant
support size, and diagonalize $R_k=\operatorname{diag}(r_k,1-r_k)$ by
output rotations, which preserve all frame identities and costs.

By \eqref{anisogap},
\[
 |\delta_k|=|2r_k-1|\le\sqrt{60\epsilon_k/11}\longrightarrow0.
\]
For an arbitrary atom put $s_i^{(k)}=\|w_i^{(k)}\|^2$.
The normalized anchor expression satisfies
\[
 g_{r_k}(a_i,b_i)\le f(s_i^{(k)})+|\delta_k|/2.
\]
This follows from $r_k(1-r_k)\le1/4$ and
$((1-r_k)a_i+r_kb_i)/2\le(1+|\delta_k|)s_i^{(k)}/4$.
As $\T^{(k)}\le g_{r_k}(a_i,b_i)$ for every atom and $f'\ge1/2$,
\[
 s_i^{(k)}\ge1-\eta_k,\qquad
 \eta_k=2\epsilon_k+|\delta_k|\longrightarrow0.
\]
Every summand below is nonnegative, and
\begin{equation}\label{energycontrol}
 \sum_i p_i^{(k)}(s_i^{(k)}-1+\eta_k)=\eta_k.
\end{equation}
For every index set $B$, summing its nonnegative terms gives
\[
 \sum_{i\in B}p_i^{(k)}s_i^{(k)}
 \le(1-\eta_k)\sum_{i\in B}p_i^{(k)}+\eta_k
 \le\sum_{i\in B}p_i^{(k)}+\eta_k.
\]
Thus atoms with vanishing probability cannot retain positive energy.

Pass to subsequences with $p_i^{(k)}\to p_i$ and
$u_i^{(k)}\to u_i$; probabilities and input directions lie in compact sets.
For $J=\{i:p_i>0\}$, the bounds
$0\le s_i^{(k)}-1+\eta_k\le\eta_k/p_i^{(k)}$ imply $s_i^{(k)}\to1$.
Pass to a further subsequence with $w_i^{(k)}\to w_i$ for $i\in J$.
For $i\notin J$, the input covariance contribution has norm $2p_i^{(k)}$,
the output covariance contribution has norm
$p_i^{(k)}s_i^{(k)}\le p_i^{(k)}+\eta_k$, and the cross-covariance contribution
has norm $\sqrt{2p_i^{(k)}}\sqrt{p_i^{(k)}s_i^{(k)}}$. All three vanish.
The surviving positive-weight atoms thus satisfy all frame identities,
have unit output norm, total energy one, and covariance $R=I_2/2$.

By Proposition~\ref{strictfinite}, this limiting frame has $\T^{(\infty)}<5/6$.
The minimum is attained because there are finitely many nonempty compact
hulls. Fix a minimizing convex combination of pair interpolators on at most
three surviving indices. Every pair used has a nonzero limiting determinant;
the same pairs are independent for all sufficiently large $k$, and their
interpolators converge. Keeping its convex coefficients gives
\[
 \limsup_k\T^{(k)}\le\T^{(\infty)}<5/6,
\]
contradicting $\T^{(k)}\to5/6$. Consequently $\Gamma<5/6$.
Set $\varepsilon_0=5/6-\Gamma>0$. This proves the uniform bound for
positive energy. If $E=0$, all $w_i$ vanish and $\T=0$, so the same
bound holds. No numerical lower bound on $\varepsilon_0$ is used.
\end{proof}

\section{Return to regression}
Let $W^\circ$ be the full-data minimum-norm minimizer and
$\rho_i=W^\circ x_i-y_i$. The normal equations imply
$\sum_i\rho_ix_i^\top=0$, and for every $B$,
\[
 L_D(W^\circ+B)=L_D^*+N^{-1}\sum_i\|Bx_i\|^2.
\]
Suppose $K=\sum_i x_ix_i^\top\succ0$.
Write $v_i=K^{-1/2}x_i$ and $\ell_i=\|v_i\|^2$.
For $\ell_i>0$, set
\[
 p_i=\ell_i/2,\qquad u_i=\sqrt{2/\ell_i}\,v_i,\qquad
 w_i=\sqrt{2/\ell_i}\,\rho_i.
\]
The omitted indices have $x_i=0$ and contribute only unavoidable loss.
The remaining atoms satisfy \eqref{frame}, with $E=\sum_{\ell_i>0}\|\rho_i\|^2\le NL_D^*$.
For $\Theta=(W-W^\circ)K^{1/2}$, the excess loss is $\|\Theta\|_F^2/N$.

Given frame weights $q$, choose original weights proportional to $q_i/\ell_i$.
The selected objective, in $\Theta$ coordinates, is a positive multiple of
$\sum_iq_i\|\Theta u_i+w_i\|^2$. If $U(q)=\sum_iq_i u_iu_i^\top\succ0$,
its unique minimizer is $\Theta=-G(q)$, where
\begin{equation}\label{mixture}
 G(q)=\left(\sum_iq_iw_iu_i^\top\right)U(q)^{-1}
 =\sum_{i<j:D_{ij}>0}\frac{q_iq_jD_{ij}}{\sum_{a<b}q_aq_bD_{ab}}G_{ij},
 \qquad D_{ij}=\det(u_i,u_j)^2.
\end{equation}
The determinant identity follows from Cauchy--Binet. By linearity of
the adjugate in dimension two, the coefficient of $q_iq_j$ in the numerator
is $M_{ij}\adj U_{ij}$, where
$U_{ij}=u_i u_i^\top+u_j u_j^\top$ and
$M_{ij}=w_i u_i^\top+w_j u_j^\top$. For an independent pair,
$M_{ij}=G_{ij}U_{ij}$, giving $D_{ij}G_{ij}$. For a dependent pair and
for diagonal terms, the adjugate annihilates the common input span.

For each eligible $S$ in the definition of $\T$, its hull is the closure
of the values $G(q)$ with $\supp q\subseteq S$ and $U(q)\succ0$
(cf.\ Zhang~\cite{zhang}). Equation~\eqref{mixture} gives one inclusion.
An independent pair realizes its single interpolator. For three pairwise
independent directions and positive coefficients
$\beta_{ij},\beta_{ik},\beta_{jk}$ of sum one, choose
\[
 q_i:q_j:q_k=
 D_{jk}/\beta_{jk}:D_{ik}/\beta_{ik}:D_{ij}/\beta_{ij}.
\]
Substitution in \eqref{mixture} gives the prescribed coefficients.
Approximate any boundary coefficient vector by positive ones.
If two directions are parallel, both valid pairs share an anchor $i$.
For a convex combination of its valid neighbors with coefficients
$\alpha_j$, take
\[
 q_i=1-t,\qquad q_j=t\frac{\alpha_j/D_{ij}}{\sum_k\alpha_k/D_{ik}},
 \qquad 0<t<1.
\]
The leaf pair has zero determinant, so \eqref{mixture} gives
$G(q)=\sum_j\alpha_jG_{ij}$ for every $0<t<1$.
All realizing and approximating fits have full rank and hence a unique
minimizer in the original coordinates. Consequently
\[
 \inf_{|\supp c|\le3}L_D(\widehat W_c)
 \le L_D^*+\T/N\le(11/6-\varepsilon_0)L_D^*.
\]

If the feature span has dimension one, let $e$ be a unit spanning vector
and write $x_i=a_ie$ and $W^\circ=\beta e^\top$. The normal equation
$\sum_i a_i\rho_i=0$ centers the vectors $\rho_i/a_i$ on nonzero $a_i$
with probabilities $a_i^2/\sum_j a_j^2$. By planar Carathéodory,
choose a convex representation of zero using at most three of them.
If its coefficients are $\lambda_i$, take $c_i\propto\lambda_i/a_i^2$.
Then $\sum_i c_i a_i\rho_i=0$ and $\sum_i c_i a_i^2>0$, so the selected
fit recovers $\beta$ on $e$. Minimum norm sets its action on $e^\perp$
to zero. If the feature span has dimension zero, choose any one observation.
All losses are independent of $W$ in that case.
In the full-rank branch, $E=0$ means that every nonzero-feature residual
vanishes, and a basis of two feature vectors recovers $W^\circ$. This covers
$L_D^*=0$ in that branch, since $E\le NL_D^*$.

For completeness, Zhang's axial example~\cite{zhang} gives the lower bound.
Take $(e_i,a+r_j)$ for the two coordinate inputs $e_i$ and the three
vertices $r_j$ of a centered equilateral triangle with $\|r_j\|^2=2/3$.
Choose $\|a\|^2=5/6$. The full-data fit has both columns equal to $a$ and
$L_D^*=2/3$. With $k=0,1,2,3$ selected points on one axis, the best squared
column deviation from $a$ is respectively $5/6,2/3,1/6,0$. For $k=0$ this
uses minimum-norm tie breaking, and for $k>0$ it follows by projecting onto
the convex hull of the chosen vertices after subtracting $a$. The excess
loss is half the sum of the two squared column deviations. The best
allocations at budget three are $(0,3)$ and $(1,2)$ up to order, both giving
excess $5/12$ and ratio $1+(5/12)/(2/3)=13/8$.

\section{Formal verification and code}
The companion source at \url{https://github.com/nZiben/weighted-selection-gap}
formalizes scalar moment selection under the two moment constraints,
its anisotropy bound \eqref{anisogap}, scalar sharpness, and the algebraic
obstruction in Proposition~\ref{strictfinite}. The planar lemma, equality
rigidity, support reduction, compactness and regression bridge are proved
in the text. The repository includes reproduction commands and the compiler
record for Lean~4.34.0 and pinned Mathlib. The audited declarations use only
\texttt{propext}, \texttt{Classical.choice} and \texttt{Quot.sound}.

\section*{Use of AI}
The proof was produced by GPT-6 Pro (GPT-6 Astra) in the ChatGPT user interface.


\begin{thebibliography}{9}
\bibitem{open}
S. Hanneke, S. Moran, A. Shlimovich and A. Yehudayoff.
\emph{Open Problem: Data Selection for Regression Tasks}.
COLT 2025, PMLR 291, pp.\ 6225--6229.
\url{https://proceedings.mlr.press/v291/hanneke25e.html}.
\bibitem{erms}
S. Hanneke, S. Moran, A. Shlimovich and A. Yehudayoff.
\emph{Data Selection for ERMs}. 2025.
\url{https://arxiv.org/abs/2504.14572}.
\bibitem{zhang}
G. Zhang.
\emph{Exact Recovery Thresholds for Weighted Data Selection in Vector-Valued Linear Regression}.
arXiv:2608.30254, version 2, 2026.
\url{https://arxiv.org/abs/2608.30254v2}.
\end{thebibliography}
\end{document}